# Touchard like polynomials and generalized Stirling numbers

G. DATTOLI †, B. GERMANO ‡, M.R. MARTINELLI
and P.E. RICCI * [1]

† Gruppo Fisica Teorica e Matematica Applicata, Unità Tecnico Scientifica
Tecnologie Fisiche Avanzate, ENEA-Centro Ricerche Frascati, C.P. 65,
Via Enrico Fermi 45, I-00044, Frascati, Roma, Italia

‡ Sapienza - Università di Roma, Dipartimento di Metodi e Modelli
Matematici per le Scienze Applicate, Via A. Scarpa, 14 - 00161 Roma, Italia

* Dipartimento di Matematica "Guido Castelnuovo", Sapienza - Università di
Roma, P.le A. Moro, 2 - 00185, Roma, Italia

The theory of Touchard polynomials is generalized using a method based on the definition of exponential operators, which extend the notion of the shift operator. The proposed technique, along with the use of the relevant operational formalism, allows the straightforward derivation of properties of this family of polynomials and their relationship to different forms of Stirling numbers.

## 1. Introduction

The use of exponential operators has been shown a powerful tool to deal with differential equations and with the properties of special functions and polynomials. In ref. [1] a thoroughly analysis of the underlying formalism has been discussed and in particular it has been shown that, under an appropriate change of variable, the action of the operator

$$\hat{E} = e^{\lambda q(x) \partial_x} \qquad (1)$$

on a given function of $x$, can be treated as an ordinary shift operator.
Just to give an example we note that in the case of the dilatation operator ($q(x) = x$) we can use the change of variables $x = e^{\vartheta}$ to prove that

$$e^{\lambda x \partial_x} f(x) = e^{\lambda \partial_\vartheta} f\left(e^{\vartheta}\right) = f\left(e^{\lambda+\vartheta}\right) = f\left(e^{\lambda} x\right). \qquad (2)$$

As a fairly direct application of the previous identity, we consider the definition of the Touchard polynomials [2] according to the following Rodrigues' type

---

[1] Corresponding author. Email: riccip@uniroma1.it

formula

$$T_n(x) = e^{-x}(x\partial_x)^n e^x.\tag{3a}$$

The relevant generating function can be obtained by multiplying both sides of eq. (3a) by $t^n/n!$ and then by summing up on the index $n$, thus getting, on account of eq. (2)

$$\sum_{n=0}^{\infty} \frac{t^n}{n!} T_n(x) = e^{-x} e^{tx\partial_x} e^x = e^{x(e^t-1)}.\tag{3b}$$

They belong therefore to the family of Sheffer polynomials, whose generating function can be written as [3]

$$\sum_{n=0}^{\infty} \frac{t^n}{n!} s_n(x) = A(t) e^{xB(t)}.\tag{4}$$

The relevant properties can therefore be studied using the technique summarized in ref. [4], where the associated multiplicative and derivative operators have been defined using the methods of monomial calculus [5]. The previous definition in terms of the operational rule (3a), allows the straightforward derivation of other relevant properties like the recurrences. From eq. (3a) we get indeed

$$(x\partial_x)\left[e^x T_n(x)\right] = (x\partial_x)^{n+1} e^x,\tag{5}$$

which can be easily manipulated to provide the recurrence

$$(x + x\partial_x)T_n(x) = T_{n+1}(x),\tag{6}$$

suggesting that the operator

$$\hat{M} = x(1 + \partial_x)\tag{7}$$

is the Touchard multiplicative operator, while the relevant derivative operator can be written as (see ref. [4])

$$\hat{P} = \ln(1 + \partial_x),\tag{8}$$

and indeed it is easily checked that

$$\ln(1 + \partial_x)T_n(x) = nT_{n-1}(x).\tag{9}$$

Before proceeding further we also remind the disentanglement identity [2]

$$e^{\hat{A}+\hat{B}} = e^{\frac{e^m-1}{m}\hat{A}} e^{\hat{B}}, \qquad \left[\hat{A}, \hat{B}\right] = -m\hat{A},\tag{10}$$

which can be exploited to get a generalization of the generating function (3b). We will indeed prove that

$$\sum_{n=0}^{\infty} \frac{t^n}{n!} T_{n+m}(x) = e^{x(e^t-1)} T_m\left(e^t x\right).\tag{11}$$

By noting indeed that

$$\hat{M}^n T_m(x) = T_{n+m}(x), \tag{12}$$

we obtain

$$\sum_{n=0}^{\infty} \frac{t^n}{n!} T_{n+m}(x) = e^{\hat{M}t} T_m(x) = e^{t(x+x\partial_x)} T_m(x). \tag{13}$$

The exponential operator in eq. (13) can now be handled by setting

$$\hat{A} = tx, \quad \hat{B} = tx\partial_x, \tag{14a}$$

and by noting that

$$\left[\hat{A}, \hat{B}\right] = -t\hat{A}. \tag{14b}$$

We finally obtain

$$\sum_{n=0}^{\infty} \frac{t^n}{n!} T_{n+m}(x) = e^{t(x+x\partial_x)} T_m(x) = e^{x(e^t-1)} T_m\left(e^t x\right). \tag{15}$$

As a by product of the previous relations, we end up with the following operational identity

$$(x + x\partial_x)^m = \sum_{r=0}^{m} \binom{m}{r} T_{m-r}(x)(x\partial_x)^r, \tag{16}$$

which can be obtained from

$$\sum_{m=0}^{\infty} \frac{t^m}{m!} (x + x\partial_x)^m = e^{x(e^t-1)} e^{tx\partial_x} = \sum_{m=0}^{\infty} \frac{t^m}{m!} \sum_{r=0}^{m} \binom{m}{r} T_{m-r}(x)(x\partial_x)^r \tag{17}$$

after equating the *t-like* powers. Before closing this section, let us remind that the following operational identity is often exploited in combinatorial analysis

$$(x\partial_x)^n = \sum_{k=0}^{n} S_2(k,n) x^k \partial_x^k, \tag{18}$$

where $S_2(k,n)$ are Stirling numbers of the second kind [6], accordingly we can obtain the following explicit definition of the Touchard polynomials

$$T_n(x) = \sum_{k=0}^{n} S_2(k,n) x^k. \tag{19}$$

In this paper we develop a systematic investigation of the properties of the generalized Touchard polynomials and of the associated Stirling numbers by the use of the operational methods we have just exploited, which will be proved to be a powerful and flexible tool of analysis.

## 2. The higher order Touchard polynomials

Before going further, we remind the formula of successive derivative of a composed functions which had different formulations in the past [7], one of these is the Hoppe's formula [8], according to which

$$\left(\frac{d}{dt}\right)^m g(f(t)) = \sum_{k=0}^{m} \frac{g^{(k)}(f(t))}{k!} A_{m,k}(t),$$

$$A_{m,k}(t) = \sum_{j=0}^{k} \binom{k}{j} (-f(t))^{k-j} \left(\frac{d}{dt}\right)^m (f(t))^j.$$
(20)

We can now use the above relation to write explicitly the Touchard polynomials, by setting indeed $x = e^\xi$ in eq. (3a), we obtain

$$\Theta_n(\xi) = e^{-e^\xi} \left(\frac{d}{d\xi}\right)^n e^{e^\xi} = \sum_{k=0}^{n} \frac{e^{k\xi}}{k!} \sum_{j=0}^{k} (-1)^{k-j} \binom{k}{j} j^n \Rightarrow$$

$$\Rightarrow T_n(x) = \sum_{k=0}^{n} \frac{x^k}{k!} \sum_{j=0}^{k} (-1)^{k-j} \binom{k}{j} j^n.$$
(21)

Which yields, after comparison with eq. (19), the well known definition of the Stirling numbers [6] as

$$S_2(k,n) = \frac{1}{k!} \sum_{j=0}^{k} (-1)^{k-j} \binom{k}{j} j^n.$$
(22)

Let us now introduce, as generalization of the Touchard polynomials, the *second order* Touchard defined as

$$T_n^{(2)}(x) = e^{-x} \left(x^2 \partial_x\right)^n e^x.$$
(23)

The analysis of the relevant properties proceeds in full analogy to what has been outlined in the introductory section for the ordinary case. We first note that in the case of exponential shift operators with $q(x) = x^2$, we obtain

$$e^{\lambda x^2 \partial_x} f(x) = f\left(\frac{x}{1-\lambda x}\right), \quad (|\lambda x| < 1).$$
(24)

as easily checked after setting $\xi = 1/x$ in the l.h.s. of eq. (24), which yields

$$e^{\lambda x^2 \partial_x} f(x) = e^{-\lambda \partial_\xi} f\left(\frac{1}{\xi}\right) = f\left(\frac{1}{\xi-\lambda}\right) = f\left(\frac{x}{1-\lambda x}\right).$$
(25)

The generating function of the second order Touchard polynomials can therefore be obtained as

$$\sum_{n=0}^{\infty} \frac{t^n}{n!} T_n^{(2)}(x) = e^{-x} e^{tx^2 \partial_x} e^x = e^{\frac{tx^2}{1-tx}}, \quad (|tx| < 1).$$
(26)

Even though the previously introduced second order Touchard polynomials do not belong to the Sheffer family, the following recursion can easily be proved

$$x^2(1+\partial_x)T_n^{(2)}(x) = T_{n+1}^{(2)}(x). \qquad (27)$$

This identity can now be exploited according to the prescription of the previous section, by noting indeed that

$$e^{\hat{A}+\hat{B}} = e^{\frac{\hat{A}}{1-\frac{m}{2}\hat{A}^{1/2}}}e^{\hat{B}}, \qquad \left[\hat{A},\hat{B}\right] = -m\hat{A}^{3/2}, \qquad (28)$$

we obtain

$$\sum_{n=0}^{\infty}\frac{t^n}{n!}T_{n+m}^{(2)}(x) = e^{t(x^2+x^2\partial_x)}T_m^{(2)}(x) = e^{\frac{tx^2}{1-tx}}T_m^{(2)}\left(\frac{x}{1-tx}\right), \quad (|tx|<1) \quad (29)$$

and

$$\left(x^2+x^2\partial_x\right)^m = \sum_{r=0}^{m}\binom{m}{r}T_{m-r}^{(2)}(x)(x^2\partial_x)^r. \qquad (30)$$

If we define, in analogy to eq. (18), the operational identity (see also ref. [9])

$$\left(x^2\partial_x\right)^n = x^n\sum_{k=0}^{n}S_2^{(2)}(k,n)x^k\partial_x^k, \qquad (31)$$

where $S_2^{(2)}(k,n)$ are a generalized form Stirling numbers, we get

$$T_n^{(2)}(x) = x^n\sum_{k=0}^{n}S_2^{(2)}(k,n)x^k. \qquad (32)$$

We can now use the same procedure as before to get the explicit expression for the Stirling numbers defined in eq. (25).
By setting indeed $\xi = 1/x$ in eq. (24), we find

$$\Theta_n^{(2)}(\xi) = e^{-\frac{1}{\xi}}\left(-\frac{d}{d\xi}\right)^n e^{\frac{1}{\xi}} = \sum_{k=0}^{n}\frac{1}{k!}\sum_{j=0}^{k}(-1)^{k-j}\binom{k}{j}\frac{(j+n)!}{(j-1)!}\frac{1}{\xi^{k+n}} \Rightarrow$$

$$\Rightarrow T_n^{(2)}(x) = x^n\sum_{k=0}^{n}\frac{x^k}{k!}\sum_{j=0}^{k}(-1)^{k-j}\binom{k}{j}\frac{(j+n-1)!}{(j-1)!}. \qquad (33)$$

The *third order* Stirling number will be therefore defined as

$$S_2^{(2)}(k,n) = \frac{1}{k!}\sum_{j=0}^{k}(-1)^{k-j}\binom{k}{j}\frac{(j+n-1)!}{(j-1)!}, \qquad (34)$$

whose properties can be studied by means of a straightforward extension of the previous procedure.

The theory of the $m$-th order Touchard polynomials can be obtained as an extension of the previous identities and the relevant properties can be summarized as reported below.

We first note that [2]

$$e^{\lambda x^m \partial_x} f(x) = f\left(\frac{x}{\sqrt[m-1]{1-(m-1)\lambda x^{m-1}}}\right), \qquad (35)$$

$$|\lambda x^{m-1}| < \frac{1}{m-1}.$$

Then by induction, we find

$$T_n^{(m)}(x) = e^{-x}(x^m \partial_x)^n e^x, \qquad (36)$$

$$(x^m + x^m \partial_x) T_n^{(m)}(x) = T_{n+1}^{(m)}(x), \qquad (37)$$

$$\sum_{n=0}^{\infty} \frac{t^n}{n!} T_{n+\ell}^{(m)}(x) = e^{\frac{x}{\sqrt[m-1]{1-(m-1)tx^{m-1}}}-x} T_\ell^{(m)}\left(\frac{x}{\sqrt[m-1]{1-(m-1)tx^{m-1}}}\right), \qquad (38)$$

$$(x^m + x^m \partial_x)^p = \sum_{r=0}^{p} \binom{p}{r} T_{p-r}^{(m)} (x^m \partial_x)^r. \qquad (39)$$

Furthermore, by setting $x = 1/\sqrt[m-1]{\xi}$ in eq. (36), we obtain

$$\Theta_n^{(m)}(\xi) = e^{-\frac{1}{m-1/\xi}}\left(-\frac{d}{d\xi}\right)^n e^{\frac{1}{m-1/\xi}} =$$

$$= \sum_{k=0}^{n} \frac{1}{k!} \sum_{j=0}^{k} (-1)^{k-j} \binom{k}{j} \frac{\Gamma\left(\frac{j}{m-1}+m+1\right)}{\Gamma\left(\frac{j}{m-1}\right)} \frac{1}{\xi^{m+\frac{k}{m-1}}} \Rightarrow \qquad (40)$$

$$\Rightarrow T_n^{(m)}(x) = \left[(m-1)x^{m-1}\right]^n \sum_{k=0}^{n} \frac{x^k}{k!} \sum_{j=0}^{k} (-1)^{k-j} \binom{k}{j} \frac{\Gamma\left(\frac{j}{m-1}+m\right)}{\Gamma\left(\frac{j}{m-1}\right)},$$

and by the use of the identity

$$T_n^{(m)}(x) = \left[(m-1)x^{m-1}\right]^n \sum_{k=0}^{n} S_2^{(m)}(k,n) x^k, \qquad (41)$$

we end up with

$$S_2^{(m)}(k,n) = \frac{1}{k!} \sum_{j=0}^{k} (-1)^{k-j} \binom{k}{j} \frac{\Gamma\left(\frac{j}{m-1}+m\right)}{\Gamma\left(\frac{j}{m-1}\right)}. \qquad (42)$$

The results we have obtained so far are quite general and offer a flexible use of the properties of this family of polynomials. In the next section we will see how they offer important links with other polynomial families.

## 3. Concluding Remarks

Albeit the tools we have provided in the previous sections allow the analysis of the generalized Touchard polynomials without any other means, it might be useful to recognize their link to other polynomials forms. We have already noted that the second order Touchard polynomials cannot be recognized as members of the Sheffer family, notwithstanding the use of the generating function of the ordinary Laguerre polynomials [10]

$$\sum_{n=0}^{\infty} t^n L_n(x) = \frac{1}{1-t} \, e^{-\frac{xt}{1-t}}, \tag{43}$$

yields the following identification of the second order Touchard in terms of Laguerre polynomials

$$T_n^{(2)}(x) = n! \left[ x^n L_n(-x) - x^{n-1} L_{n-1}(-x) \right]. \tag{44}$$

Identifications of this type becomes more difficult with increasing order, however the link with the Bessel polynomials $b_n(x)$ firstly introduced by Krall and Frinck [11] is noticeable and worth to be stressed.
Their properties are well known and well documented, they have been studied within the context of Sheffer series (see ref. [6]) and in this section we will see how they can be framed within the present formalism.
We consider indeed the polynomials defined through the following Rodriguez type relation

$$\delta_n(x) = e^{-1/x} \left( x^3 \partial_x \right)^n e^{1/x}. \tag{45}$$

The use of the operational methods we have described in the previous sections allow to cast the above polynomial family in the form

$$\sum_{n=0}^{\infty} \frac{(-t)^n}{n!} \delta_n(-x) = e^{\frac{1}{x}\left(1-\sqrt{1-2tx}\right)}, \tag{46}$$

thus recognizing them with the Bessel polynomials according to the following identification

$$\delta_n(-x) = x^{2n} y_{n-1} \left( \frac{1}{x} \right). \tag{47}$$

In a forthcoming investigation we will reconsider the present result and their usefulness within the context of different disciplines, including combinatorial analysis and quantum optics.

Author's addresses

Giuseppe Dattoli
Unità Tecnico Scientifica Tecnologie Fisiche Avanzate
ENEA – Centro Ricerche Frascati – C.P. 65
Via E. Fermi, 45
00044 - Frascati - Roma (Italia)
e-mail: dattoli@frascati.enea.it

Bruna Germano
Dipartimento di Metodi e Modelli Matematici per le Scienze Applicate
Università degli Studi di Roma "La Sapienza"
Via A. Scarpa, 14
00161 – Roma (Italia)
e-mail: germano@dmmm.uniroma1.it

Maria Renata Martinelli
Dipartimento di Metodi e Modelli Matematici per le Scienze Applicate
Università degli Studi di Roma "La Sapienza"
Via A. Scarpa, 14
00161 – Roma (Italia)
e-mail: martinelli@dmmm.uniroma1.it

Paolo E. Ricci
Dipartimento di Matematica "Guido Castelnuovo"
Sapienza – Università di Roma
P.le A. Moro, 2
00185 – Roma (Italia)
e-mail: riccip@uniroma1.it